\documentclass[12pt]{article}

\makeatletter
\renewcommand{\@oddfoot}{\hfill \thepage}
\makeatother

\usepackage[T2A]{fontenc}
\usepackage[cp1251]{inputenc}
\usepackage[russian, english]{babel}
\usepackage{mathrsfs}
\usepackage{amssymb}
\usepackage{amsmath,amsthm,amsfonts}
\usepackage{cite}
\usepackage[dvipdf]{graphicx}
\usepackage{fancyhdr}

\begin{document}

\begin{center}
{\bf LINEAR FORMS OF THE TELEGRAPH RANDOM PROCESSES\\ DRIVEN BY PARTIAL DIFFERENTIAL EQUATIONS}
\end{center}

\begin{center}
Alexander D. KOLESNIK\\
Institute of Mathematics and Computer Science\\
Academy Street 5, Kishinev 2028, Moldova\\
E-Mail: kolesnik@math.md
\end{center}

\vskip 0.2cm

\begin{abstract}
Consider $n$ independent Goldstein-Kac telegraph processes $X_1(t), \dots ,X_n(t), \; n\ge 2, \; t\ge 0,$ on the real line $\Bbb R$.
Each the process $X_k(t), \; k=1,\dots,n,$ describes a stochastic motion at constant finite speed $c_k>0$  of a particle 
that, at the initial time instant $t=0$, starts from some initial point $x_k^0=X_k(0)\in\Bbb R$
and whose evolution is controlled by a homogeneous Poisson process $N_k(t)$ of rate $\lambda_k>0$. The governing Poisson processes
$N_k(t), \; k=1,\dots,n,$ are supposed to be independent as well. Consider the linear form of the processes
$X_1(t), \dots ,X_n(t), \; n\ge 2,$ defined by
$$L(t) = \sum_{k=1}^n a_k X_k(t) ,$$
where $a_k, \; k=1,\dots,n,$ are arbitrary real non-zero constant coefficients.
We obtain a hyperbolic system of first-order partial differential equations for the joint probability
densities of the process $L(t)$ and of the directions of motions at arbitrary time $t>0$. From this system
we derive a partial differential equation of order $2^n$ for the transition density of $L(t)$ in the form of a determinant of a block matrix whose elements 
are the differential operators with constant coefficients. The weak convergence of $L(t)$ to a homogeneous Wiener process, under Kac's scaling conditions, 
is proved. Initial-value problems for the transition densities of the sum and difference $S^{\pm}(t)=X_1(t) \pm X_2(t)$ of two independent telegraph 
processes with arbitrary parameters, are also posed. 
\end{abstract}

\vskip 0.1cm


\vskip 0.2cm


\section{Introduction}

\numberwithin{equation}{section}

The classical Goldstein-Kac telegraph process, first introduced in the works \cite{gold} and \cite{kac}, describes the stochastic motion of a particle that moves at constant finite speed on the real line $\Bbb R$ and alternates two possible directions of motion at random Poisson time instants. The main properties of this process and its numerous generalizations, as well as some their applications, have been studied in a series of works \cite{bart1,bart2,br,cane1,cane2,cres1,cres2,foong1,foong2}, \cite{iacus1,iacus2,kab}, \cite{kap,kis}, \cite{kol3,kol4}, \cite{kol6,loprat,mas1,mas2,pin,rat1,rat2,sta}, \cite{turb}. An introduction to the contemporary theory of the telegraph processes and their applications in financial modelling can be found in the recently published book \cite{kol6}.

In all of the above mentioned works the main subject of interest was a single particle whose evolution develops at some finite speed on $\Bbb R$. On the other hand, studying a system of several telegraph particles is of a special importance from the point of view of describing various kinds of interactions (for the interpretation of different interactions arising in physics, chemistry, biology, financial markets, see \cite[p. 1173]{kol2}). However, to the best of the author's knowledge, there exist only a few works where a system of several telegraph processes is considered. A closed-form expression for the probability distribution function of the Euclidean distance between two independent telegraph processes with arbitrary parameters was derived in \cite{kol2}. The explicit probability distribution for the sum of two independent telegraph processes with the same parameters, both starting from the origin of $\Bbb R$, was obtained in \cite{kol1}. 

In this article we continue studying the problems devoted to the evolutions of several particles moving randomly with finite speed. The subject of our interests is the general linear form 
\begin{equation}\label{form1}
L(t) = \sum_{k=1}^n a_k X_k(t), \qquad  a_k\in\Bbb R, \quad a_k\neq 0, \;\; k=1,\dots,n,  \quad t\ge 0,
\end{equation}
of $n$ independent telegraph processes $X_1(t), \dots ,X_n(t), \; n\ge 2, \; t\ge 0$, that, at the initial time instant $t=0$, simultaneously start from arbitrary initial points $x_k^0\in\Bbb R, \; k=1,\dots, n$. Clearly, linear form (\ref{form1}) is a generalization of the sum of two telegraph processes with the same parameters studied in \cite{kol1}. 

The process $L(t)$ defined by (\ref{form1}) is of a great interest since it has a quite natural interpretation as a weighted linear combination of $n$ independent stochastic sources of the telegraph type. From this point of view, the knowledge of the probabilistic characteristics of process $L(t)$ enables to study various important problems, such as attaining a given level by the process or being in a fixed interval. Besides purely mathematical interest, this interpretation justifies the strong motivation for studying such type of processes. 

The paper is organized as follows. In Section 2 we recall some basic properties of the telegraph process that 
we will substantially be relying on. In Section 3 we describe the structure of the distribution of linear form (\ref{form1}) and derive a hyperbolic system of 
first-order partial differential equations for the joint densities of process (\ref{form1}) and of the directions of motion 
at arbitrary time $t>0$. From this system, a hyperbolic partial differential equation of order $2^n$ with constant coefficients for the transition probability density of $L(t)$ is given in Section 4 in the form of a determinant having a nice block structure whose elements are the commuting first-order differential operators of the system.
In Section 5 we prove that, under Kac's scaling conditions, the process $L(t)$ weakly converges to a homogeneous Brownian motion on $\Bbb R$. Initial-value problems for the sum and difference $S^{\pm}(t)=X_1(t) \pm X_2(t)$ of two independent telegraph processes with arbitrary parameters, are posed in Section 6.

\section{Some Basic Properties of the Telegraph Process}

\numberwithin{equation}{section}

The Goldstein-Kac telegraph process describes the stochastic motion of a particle that, at the initial time instant $t=0$, starts from the origin $x=0$ of the real
line $\Bbb R$ and moves with some finite constant speed $c$. The initial direction of the motion (positive or negative) is taken at random  
with equal probabilities 1/2. The motion is controlled by a homogeneous Poisson process of rate $\lambda>0$ as follows. At every 
Poisson event occurence, the particle instantaneously takes on the opposite direction and keeps moving with the same speed $c$ until
the next Poisson event occurs, then it takes on the opposite direction again, and so on.
This random motion was first studied in the works by Goldstein \cite{gold} and Kac \cite{kac}.

Let $X(t)$ denote the particle's position on $\Bbb R$ at arbitrary time instant $t>0$. Since the speed $c$ is finite, then 
at time $t>0$ the distribution $\text{Pr}\{ X(t)\in dx \}$ is concentrated in the finite interval $[-ct, ct]$ which is
the support of this distribution. The density $f(x,t), \; x\in\Bbb R, \; t\ge 0,$ of the distribution $\text{Pr}\{
X(t)\in dx \}$ has the structure $$f(x, t) = f_{s}(x, t) + f_{ac}(x, t),$$
where $f_{s}(x, t)$ and $f_{ac}(x, t)$ are the densities of the singular (with respect to the Lebesgue measure on the line) and of
the absolutely continuous components of the distribution of $X(t)$, respectively.

The singular component is obviously concentrated at two terminal points $\pm ct$ of the interval
$[-ct, ct]$ and corresponds to the case when no one Poisson event occurs until the moment $t$ and, therefore, the particle does not
change its initial direction. Therefore, the probability of being at time $t>0$ at the terminal points $\pm ct$ is
\begin{equation}\label{prop1}
\text{Pr}\left\{ X(t) = ct \right\} = \text{Pr}\left\{ X(t)
= - ct \right\} = \frac{1}{2} \; e^{-\lambda t} .
\end{equation}

The absolutely continuous component of the distribution of $X(t)$ is concentrated in the open interval $(-ct, ct)$ and corresponds
to the case when at least one Poisson event occurs before time instant $t$ and, therefore, the particle changes its initial direction.
The probability of this event is
\begin{equation}\label{prop2}
\text{Pr}\left\{ X(t) \in (-ct, ct) \right\} = 1 - e^{-\lambda t}.
\end{equation}

The principal result by Goldstein \cite{gold} and Kac \cite{kac}
states that the density $f = f(x,t), \; x\in [-ct, ct], \; t>0,$
of the distribution of $X(t)$ satisfies the hyperbolic partial differential equation
\begin{equation}\label{prop3}
\frac{\partial^2 f}{\partial t^2} + 2\lambda
\frac{\partial f}{\partial t} - c^2 \frac{\partial^2 f}{\partial x^2} = 0,
\end{equation}
(which is referred to as the {\it telegraph} or {\it damped wave}
equation) and can be found by solving (\ref{prop3}) with the initial conditions
\begin{equation}\label{iprop3}
f(x,t)\vert_{t=0} = \delta(x), \qquad
\left.\frac{\partial f(x,t)}{\partial t}\right\vert_{t=0} = 0,
\end{equation}
where $\delta(x)$ is the Dirac delta-function. This means that the transition density $f(x,t)$ of the process $X(t)$ is the
fundamental solution (i.e. the Green's function) of the telegraph equation (\ref{prop3}).

The explicit form of the density $f(x,t)$ is given by the formula (see, for instance, \cite[Section 2.5]{kol6} or \cite[Section 0.4]{pin}):
\begin{equation}\label{prop4}
\aligned
f(x,t) & = \frac{e^{-\lambda t}}{2} \left[ \delta(ct-x) +
\delta(ct+x) \right]\\
& \quad + \frac{e^{-\lambda t}}{2c} \left[ \lambda I_0\left(
\frac{\lambda}{c} \sqrt{c^2t^2-x^2} \right) +
\frac{\partial}{\partial t} I_0\left( \frac{\lambda}{c}
\sqrt{c^2t^2-x^2} \right) \right] \Theta(ct-\vert x\vert),
\endaligned
\end{equation}
where $\Theta(x)$ is the Heaviside step function
\begin{equation}\label{pprop4}
\Theta(x) = \left\{ \aligned 1, \qquad  & \text{if} \; x>0,\\
                               0, \qquad & \text{if} \; x\le 0.
\endaligned \right.
\end{equation}
Here $I_0(z)$ is the modified Bessel function of order zero (that is, the Bessel function with imaginary argument) given by the formula
\begin{equation}\label{pprop5}
I_0(z) = \sum_{k=0}^{\infty} \frac{1}{(k!)^2}
\left( \frac{z}{2} \right)^{2k} .
\end{equation}
The first term in (\ref{prop4})
\begin{equation}\label{prop5}
f^s(x,t) = \frac{e^{-\lambda t}}{2} \left[ \delta(ct-x) +
\delta(ct+x) \right]
\end{equation}
represents the density (as a generalized function) of the singular part of the distribution
of $X(t)$ concentrated at two terminal points $\pm ct$ of the interval $[-ct, ct]$, while the second term in (\ref{prop4})
\begin{equation}\label{prop6}
f^{ac}(x,t) = \frac{\lambda e^{-\lambda t}}{2c}
\left[ \lambda I_0\left( \frac{\lambda}{c} \sqrt{c^2t^2-x^2} \right) +
\frac{\partial}{\partial t} I_0\left( \frac{\lambda}{c}
\sqrt{c^2t^2-x^2} \right) \right] \Theta(ct-\vert x\vert),
\end{equation}
is the density of the absolutely continuous part of the distribution concentrated in the open interval $(-ct, ct)$.

Denote by $X^{x^0}(t)$ the telegraph process starting from arbitrary initial point $x^0\in\Bbb R$. It is clear that
the transition density of $X^{x^0}(t)$ emerges from (\ref{prop4}) by the replacement $x \mapsto x-x^0$ and has the form
\begin{equation}\label{prop7}
\aligned
f^{x^0}(x,t) & = \frac{e^{-\lambda t}}{2} \left[ \delta(ct-(x-x^0)) +
\delta(ct+(x-x^0)) \right] \\
& + \frac{e^{-\lambda t}}{2c}
\left[ \lambda I_0\left( \frac{\lambda}{c} \sqrt{c^2t^2-(x-x^0)^2} \right)
+ \frac{\partial}{\partial t} I_0\left( \frac{\lambda}{c} \sqrt{c^2t^2-(x-x^0)^2} \right)
\right] \Theta(ct-\vert x-x^0\vert).
\endaligned
\end{equation}
The support of the distribution of $X^{x^0}(t)$ is the close interval $[x^0-ct, x^0+ct]$. The first term in (\ref{prop7})
\begin{equation}\label{prop8}
f^{x^0}_s(x,t) = \frac{e^{-\lambda t}}{2} \left[ \delta(ct-(x-x^0)) +
\delta(ct+(x-x^0)) \right]
\end{equation}
is the singular part of the density concentrated at the two terminal points $x^0\pm ct$ of the support, while the second term
\begin{equation}\label{prop9}
\aligned
f^{x^0}_{ac}(x,t) & = \frac{e^{-\lambda t}}{2c}
\left[ \lambda I_0\left( \frac{\lambda}{c} \sqrt{c^2t^2-(x-x^0)^2}
\right)\right.\\
& \qquad\qquad\qquad\qquad \left. + \frac{\partial}{\partial t} I_0\left( \frac{\lambda}{c} \sqrt{c^2t^2-(x-x^0)^2} \right)
\right] \Theta(ct-\vert x-x^0\vert), \endaligned
\end{equation}
is the density of the absolutely continuous part of the distribution of $X^{x^0}(t)$ concentrated in the open interval $(x^0-ct, x^0+ct)$.

The characteristic function of the telegraph process starting from the origin $x=0$ with density (\ref{prop4}) is given by the formula 
(see \cite[Section 2.4]{kol6}):
\begin{equation}\label{prop10}
\aligned
H(\alpha, t) = e^{-\lambda t} & \left\{ \left[ \cosh\left( t\sqrt{\lambda^2-c^2\alpha^2} \right) +
\frac{\lambda}{\sqrt{\lambda^2-c^2\alpha^2}} \; \sinh\left(
t\sqrt{\lambda^2-c^2\alpha^2} \right) \right] \bold 1_{\left\{\vert\alpha\vert\le\frac{\lambda}{c}\right\}} \right. \\
& + \left. \left[ \cos\left( t\sqrt{c^2\alpha^2 - \lambda^2}
\right) + \frac{\lambda}{\sqrt{c^2\alpha^2-\lambda^2}} \;
\sin\left( t\sqrt{c^2\alpha^2-\lambda^2} \right)  \right] \bold 1_{\left\{\vert\alpha\vert > \frac{\lambda}{c}\right\}} \right\} ,
\endaligned
\end{equation}
where $\bold 1_{ \{ \cdot \}}$ is the indicator function, $\alpha\in\Bbb R, \; t\ge 0$.
Clearly, if the process starts from some arbitrary point $x^0\in\Bbb R$ and has the transition density (\ref{prop7}),
then its characteristic function $H^{x^0}(\alpha, t)$ expresses through (\ref{prop10}) as follows:
\begin{equation}\label{prop11}
H^{x^0}(\alpha, t) = e^{i\alpha x^0} H(\alpha, t) , \qquad \alpha\in\Bbb R, \quad t\ge 0.
\end{equation}

\section{Structure of Distribution and System of Equations} 

\numberwithin{equation}{section}

Let $X_1^{x_1^0}(t),\dots, X_n^{x_n^0}(t), \; n\ge 2, \; t\ge 0,$ be independent Goldstein-Kac telegraph processes on the real line $\Bbb R$ that, at the initial time instant $t=0$, simultaneously start from the initial points $x_1^0, \dots, x_n^0\in\Bbb R$, respectively. For the sake of simplicity, we omit thereafter the upper indices by identifying $X_k(t) \equiv X_k^{x_k^0}(t), \; k=1,\dots, n,$ bearing in mind, however, the fact that the process $X_k(t)$ starts from the initial point $x_k^0$. 
Each process $X_k(t), \; k=1,\dots,n,$ has some constant finite speed $c_k>0$  and is controlled by a homogeneous Poisson process $N_k(t)$ of rate $\lambda_k>0$, as described above.  
All these Poisson processes $N_k(t), \; k=1,\dots,n,$ are supposed to be independent as well. Consider the linear form of the processes
$X_1(t), \dots ,X_n(t), \; n\ge 2,$ defined by the equality
\begin{equation}\label{syst1}
L(t) = \sum_{k=1}^n a_k X_k(t) , \qquad a_k\in\Bbb R, \quad a_k\neq 0, \;\; k=1,\dots,n,  \quad t\ge 0, 
\end{equation}
where $a_k, \; k=1,\dots,n,$ are arbitrary real non-zero constant coefficients. 

To describe the structure of the distribution of $L(t)$, consider the following partition of the set of indices: 
$$\aligned 
I^+ & = \{ i_1,\dots,i_k \} \; \text{such that} \; a_{i_s}>0 \; \text{for all} \; i_s\in I^+, \; 1\le s\le k, \\
I^- & = \{ i_1,\dots,i_m \} \; \text{such that} \; a_{i_l}<0 \; \text{for all} \; i_l\in I^-, \; 1\le l\le m,  \qquad k+m=n.
\endaligned$$
The support of the distribution $\Phi(x,t)=\text{Pr}\{ L(t)<x \}$ of process 
$L(t)$ is the close interval depending on the coefficients $a_k$, speeds $c_k$ and start points $x_k^0$ and having the form: 
\begin{equation}\label{supp1}
\text{supp} \; L(t) = \left[ \sum_{k=1}^n a_kx_k^0 - t\biggl( \sum_{i_s\in I^+} a_{i_s}c_{i_s} - \sum_{i_l\in I^-} a_{i_l}c_{i_l} \biggr) ,  \;\; 
\sum_{k=1}^n a_kx_k^0 + t\biggl( \sum_{i_s\in I^+} a_{i_s}c_{i_s} - \sum_{i_l\in I^-} a_{i_l}c_{i_l} \biggr)  \right] .
\end{equation}
In particular, if all $a_k>0, \; k=1,\dots,n$, then the set $I^-$ is empty and, therefore, support (\ref{supp1}) takes the form
\begin{equation}\label{supp2}
\text{supp} \; L(t) = \left[ \sum_{k=1}^n a_k(x_k^0-c_kt), \; \sum_{k=1}^n a_k(x_k^0+c_kt) \right] .
\end{equation}

At arbitrary time instant $t>0$, the distribution $\Phi(x,t)$ contains the singular and absolutely continuous components. 
The singular part of the distribution corresponds to the case, when no one Poisson event (of any Poisson process $N_k(t), \; k=1,\dots,n,$) occurs by time instant $t$. 
It is concentrated in the finite point set $M_s = \{ q_1, \dots, q_{2^n}\}\subset\text{supp} \; L(t)$ that contains $2^n$ singularity points (each $q_j$ is counted according to its multiplicity):
\begin{equation}\label{sing1}
q_j = \sum_{k=1}^n a_kx_k^0 + t\sum_{k=1}^n a_k i_k^j c_k , \qquad j=1,\dots, 2^n,
\end{equation}
where $i_k^j=\pm 1, \; k=1,\dots,n,$ are the elements of the ordered sequence $\boldsymbol\sigma_j=\{ i_1^j,\dots, i_n^j \}, \; j=1,\dots, 2^n,$ of length $n$. The sign of each $i_k^j$, (that is, either $+1$ or $-1$), is determined by the initial direction (either positive or negative, respectively) taken by the telegraph process $X_k(t)$. Emphasize that some $q_j$ may coincide in dependence on the particular values of the start points $x_k^0$, coefficients $a_k$ and speeds $c_k$.

Note that both the terminal points of support (\ref{supp1}) are singular and, therefore, they belong to $M_s$, that is, 
$$\sum_{k=1}^n a_kx_k^0 \pm t\biggl( \sum_{i_s\in I^+} a_{i_s}c_{i_s} - \sum_{i_l\in I^-} a_{i_l}c_{i_l} \biggr) \in M_s .$$
Other singular points are the interior points of support (\ref{supp1}). It is easy to see that the probability of being at arbitrary singularity point $q_j$ 
(taking into account its multiplicity) at time instant $t$ is
\begin{equation}\label{sing2}
\text{Pr} \left\{ L(t) = q_j \right\} = \frac{e^{-\boldsymbol\lambda t}}{2^n} , \qquad j=1, \dots, 2^n, 
\end{equation}
where   
\begin{equation}\label{lam}
\boldsymbol\lambda = \sum_{k=1}^n \lambda_k.
\end{equation}
From (\ref{sing2}) it obviously follows that, for arbitrary $t>0$,  
\begin{equation}\label{sing3}
\text{Pr} \left\{ L(t) \in M_s \right\} = e^{-\boldsymbol\lambda t} .  
\end{equation}

If at least one Poisson event occurs by time instant $t$, then the process $L(t)$ is located in the set $M_{ac} = \text{supp} \; L(t) - M_s,$
which is the support of the absolutely continuous part of the distribution and the probability of being in this set at time instant $t>0$ is: 
\begin{equation}\label{cont}
\text{Pr} \left\{ L(t) \in M_{ac} \right\} = 1 - e^{-\boldsymbol\lambda t} .  
\end{equation}

Define now the two-state direction processes $D_1(t),\dots, D_n(t), \; n\ge 2, \; t>0,$
where $D_k(t), \; k=1,\dots,n,$ denotes the direction of the telegraph process $X_k(t)$ at instant $t>0$.
This means that $D_k(t)=+1,$ if at instant $t$ the process $X_k(t)$ is developing in the positive direction and $D_k(t)=-1$ otherwise.

Introduce the joint probability densities of the process $L(t)$ and of the set of directions
$ \left\{ D_1(t),\dots,D_n(t) \right\} $ at arbitrary time instant $t>0$ by the relation
\begin{equation}\label{syst2}
f_{\boldsymbol\sigma}(x,t) \; dx \equiv f_{ \{i_1,\dots,i_n\} }(x,t) \; dx = \text{Pr} \{  x<L(t)<x+dx, \; D_1(t)=i_1, \dots, D_n(t)=i_n \} .
\end{equation}
The set of functions (\ref{syst2}) contains $2^n$ densities indexed by all the ordered sequences of the form $ \boldsymbol\sigma= \{ i_1,\dots,i_n \} $
of length $n$ whose elements $i_k, \; k=1,\dots,n$, are either $+1$ or $-1$.

Our first result is given by the following theorem.

\bigskip

{\bf Theorem 1.} {\it The joint probability densities} (\ref{syst2}) {\it satisfy the following  
hyperbolic system of $2^n$ first-order partial differential equations with constant coefficients:}
\begin{equation}\label{syst3}
\frac{\partial f_{\boldsymbol\sigma}(x,t)}{\partial t} =  - \bold c_{\boldsymbol\sigma}
\frac{\partial f_{\boldsymbol\sigma}(x,t)}{\partial x} - \boldsymbol\lambda f_{\boldsymbol\sigma}(x,t)
+ \sum_{k=1}^n \lambda_k \; f_{\boldsymbol{\bar\sigma}^{(k)}}(x,t) \; .
\end{equation}
{\it where} 
$$\aligned 
\boldsymbol\sigma & = \{ i_1,\dots,i_{k-1},i_k,i_{k+1},\dots,i_n \} , \\ 
\boldsymbol{\bar\sigma}^{(k)} & = \{ i_1,\dots,i_{k-1},-i_k,i_{k+1},\dots,i_n \}, \endaligned$$  

\begin{equation}\label{coeffC}
\bold c_{\boldsymbol\sigma} \equiv \bold c_{ \{i_1,\dots,i_n\} } = \sum_{k=1}^n  a_k i_k c_k ,
\end{equation}
{\it and $\boldsymbol\lambda$ is given by} (\ref{lam}).

\vskip 0.2cm

{\it Proof.} Let $\Delta t>0$ be some time increment. Let $N_k(t,t+\Delta t), \; k=1,\dots,n,$ denote the number of the events of $k$-th Poisson
process $N_k(t)$ that have occurred in the time interval $(t,t+\Delta t)$. Then, according to the total probability formula, we have:
$$\aligned
\text{Pr} & \{  L(t+\Delta t)<x, \; D_1(t+\Delta t)=i_1, \dots, D_n(t+\Delta t)=i_n \} \\
& = \prod_{k=1}^n (1- \lambda_k \Delta t) \; \text{Pr} \left\{  L(t) + \Delta t \sum_{k=1}^n a_k i_k c_k < x,
\; D_1(t)=i_1, \dots, D_n(t)=i_n \right\} \\
& + \sum_{k=1}^n \lambda_k \Delta t \; \prod_{\substack{j=1\\ j\neq k}}^n \left( 1 - \lambda_j \Delta t \right)
\frac{1}{\Delta t} \int\limits_t^{t+\Delta t} \text{Pr} \left\{  L(t) + a_k c_k(-i_k(\tau_k-t)+i_k(t+\Delta t-\tau_k)) < x, \right.\\
& \left. D_1(t)=i_1, \dots, D_{k-1}(t)=i_{k-1}, D_k(t)=-i_k, D_{k+1}=i_{k+1}, \dots, D_n(t)=i_n \right\} \; d\tau_k + o(\Delta t) \\
& = \prod_{k=1}^n (1- \lambda_k \Delta t) \; \text{Pr} \left\{  L(t) < x - \bold c_{\boldsymbol\sigma}\Delta t , \; D_1(t)=i_1, \dots, D_n(t)=i_n \right\} \\
& + \sum_{k=1}^n \lambda_k \; \prod_{\substack{j=1\\ j\neq k}}^n \left( 1 - \lambda_j \Delta t \right)
\int\limits_t^{t+\Delta t} \text{Pr} \left\{  L(t) < x - a_k i_k c_k(2(t-\tau_k)+\Delta t), \right.\\
& \left. D_1(t)=i_1, \dots, D_{k-1}(t)=i_{k-1}, D_k(t)=-i_k, D_{k+1}=i_{k+1}, \dots, D_n(t)=i_n \right\} \; d\tau_k + o(\Delta t) .
\endaligned$$
The first term on the right-hand side of this expression is related to the case when no one Poisson event has occurred in the time interval
$(t,t+\Delta t)$, that is, if $\sum\limits_{k=1}^n N_k(t,t+\Delta t) = 0$. The second (integral) term concerns the case when a single Poisson event has occurred
in this interval, that is, if $\sum\limits_{k=1}^n N_k(t,t+\Delta t) = 1$.  Finally, the term $o(\Delta t)$ is related to the case when more that one Poisson
events have occurred in the interval $(t,t+\Delta t)$, that is, if $\sum\limits_{k=1}^n N_k(t,t+\Delta t) \ge 2$ (one can easily check that all such
probabilities have the order $o(\Delta t)$).

Since the probability is a continuous function, then, according to the mean-value theorem of classical analysis, for any $k$ there exists
a time instant $\tau_k^*\in (t,t+\Delta t)$, such that
$$\aligned
\text{Pr} & \{  L(t+\Delta t)<x, \; D_1(t+\Delta t)=i_1, \dots, D_n(t+\Delta t)=i_n \} \\
& = \prod_{k=1}^n (1 - \lambda_k \Delta t) \; \text{Pr} \left\{  L(t) < x - \bold c_{\boldsymbol\sigma}\Delta t, \; D_1(t)=i_1, \dots, D_n(t)=i_n \right\} \\
& + \Delta t \sum_{k=1}^n \lambda_k \; \prod_{\substack{j=1\\ j\neq k}}^n \left( 1 - \lambda_j \Delta t \right)
\; \text{Pr} \left\{  L(t) < x - a_k i_k c_k(2(t-\tau_k^*)+\Delta t), \right.\\
& \left. D_1(t)=i_1, \dots, D_{k-1}(t)=i_{k-1}, D_k(t)=-i_k, D_{k+1}=i_{k+1}, \dots, D_n(t)=i_n \right\} + o(\Delta t) .
\endaligned$$
In view of the asymptotic formulas
$$ \prod_{j=1}^n \left( 1 - \lambda_j \Delta t \right) = 1 - \boldsymbol\lambda\Delta t + o(\Delta t), \qquad
\Delta t \prod_{\substack{j=1\\ j\neq k}}^n \left( 1 - \lambda_j \Delta t \right) = \Delta t + o(\Delta t),$$
the latter relation can be rewritten as follows:
$$\aligned
\text{Pr} & \{  L(t+\Delta t)<x, \; D_1(t+\Delta t)=i_1, \dots, D_n(t+\Delta t)=i_n \} \\
& = \text{Pr} \left\{  L(t) < x - \bold c_{\boldsymbol\sigma}\Delta t , \; D_1(t)=i_1, \dots, D_n(t)=i_n \right\} \\
& \quad - \boldsymbol\lambda\Delta t \; \text{Pr} \left\{ L(t) < x - \bold c_{\boldsymbol\sigma}\Delta t , \; D_1(t)=i_1, \dots, D_n(t)=i_n \right\} \\
& \quad + \Delta t \sum_{k=1}^n \lambda_k \; \text{Pr} \left\{  L(t) < x - a_k i_k c_k(2(t-\tau_k^*)+\Delta t), \right.\\
& \quad \left. D_1(t)=i_1, \dots, D_{k-1}(t)=i_{k-1}, D_k(t)=-i_k, D_{k+1}=i_{k+1}, \dots, D_n(t)=i_n \right\} + o(\Delta t) .
\endaligned$$
In terms of densities (\ref{syst2}) this equality can be represented in the form:
$$\aligned
\int\limits_{-\infty}^x f_{\boldsymbol\sigma} (\xi,t+\Delta t) \; d\xi & = \int\limits_{-\infty}^{x-\bold c_{\boldsymbol\sigma}\Delta t}
f_{\boldsymbol\sigma}(\xi,t) \; d\xi
- \boldsymbol\lambda\Delta t \int\limits_{-\infty}^{x-\bold c_{\boldsymbol\sigma}\Delta t} f_{\boldsymbol\sigma}(\xi,t) \; d\xi \\
& + \Delta t \sum_{k=1}^n \lambda_k \int\limits_{-\infty}^{x - a_k i_k c_k(2(t-\tau_k^*)+\Delta t)} f_{\boldsymbol{\bar\sigma}^{(k)}}(\xi,t)
\; d\xi + o(\Delta t) .
\endaligned$$
This can be rewritten as follows:
$$\aligned
\int\limits_{-\infty}^x & \bigl[ f_{\boldsymbol\sigma} (\xi,t+\Delta t) - f_{\boldsymbol\sigma} (\xi,t) \bigr] d\xi
= - \left[\int\limits_{-\infty}^x f_{\boldsymbol\sigma}(\xi,t) - \int\limits_{-\infty}^{x-\bold c_{\boldsymbol\sigma}\Delta t}
f_{\boldsymbol\sigma}(\xi,t) \; d\xi \right] \\
& - \boldsymbol\lambda\Delta t \int\limits_{-\infty}^{x-\bold c_{\boldsymbol\sigma}\Delta t} f_{\boldsymbol\sigma}(\xi,t) \; d\xi
+ \Delta t \sum_{k=1}^n \lambda_k \int\limits_{-\infty}^{x - a_k i_k c_k(2(t-\tau_k^*)+\Delta t)} f_{\boldsymbol{\bar\sigma}^{(k)}}
(\xi,t) \; d\xi + o(\Delta t) .
\endaligned$$
Dividing this equality by $\Delta t$, we can represent it in the form:
$$\aligned
\int\limits_{-\infty}^x & \frac{1}{\Delta t} \bigl[ f_{\boldsymbol\sigma}(\xi,t+\Delta t) - f_{\boldsymbol\sigma}(\xi,t) \bigr] d\xi
= - \bold c_{\boldsymbol\sigma} \left\{ \frac{1}{\bold c_{\boldsymbol\sigma}\Delta t} \left[\int\limits_{-\infty}^x f_{\boldsymbol\sigma}(\xi,t) -
\int\limits_{-\infty}^{x-\bold c_{\boldsymbol\sigma}\Delta t} f_{\boldsymbol\sigma}(\xi,t) \; d\xi \right]\right\} \\
& \qquad\qquad - \boldsymbol\lambda \int\limits_{-\infty}^{x-\bold c_{\boldsymbol\sigma}\Delta t} f_{\boldsymbol\sigma}(\xi,t) \; d\xi
+ \sum_{k=1}^n \lambda_k \int\limits_{-\infty}^{x - a_k i_k c_k(2(t-\tau_k^*)+\Delta t)} f_{\boldsymbol{\bar\sigma}^{(k)}}(\xi,t) \; d\xi
+ \frac{o(\Delta t)}{\Delta t} .
\endaligned$$
Passing now to the limit, as $\Delta t\to 0$, and taking into account that $\tau_k^*\to t$ in this case, we obtain
$$\int\limits_{-\infty}^x \frac{\partial f_{\boldsymbol\sigma}(\xi,t)}{\partial t} \; d\xi
= - \bold c_{\boldsymbol\sigma} f_{\boldsymbol\sigma}(x,t) - \boldsymbol\lambda\int\limits_{-\infty}^x f_{\boldsymbol\sigma}(\xi,t) \; d\xi
+ \sum_{k=1}^n \lambda_k \int\limits_{-\infty}^x f_{\boldsymbol{\bar\sigma}^{(k)}}(\xi,t) \; d\xi . $$
Differentiating this equality in $x$, we finally arrive at (\ref{syst3}).

Since the principal part of system (\ref{syst3}) is strictly hyperbolic, then system (\ref{syst3}) itself is hyperbolic. 
The theorem is thus completely proved. $\square$

\bigskip

{\it Remark 1.} Note that system (\ref{syst3}) consists of $2^n$ first-order partial differential equations, however the equation for each density
$f_{\boldsymbol\sigma}(x,t)$ contains (besides this function itself) only $n$ other densities $f_{\boldsymbol{\bar\sigma}^{(k)}}(x,t), \; k=1,\dots,n$.
This means that each density $f_{\boldsymbol\sigma}(x,t)$ indexed by some ordered sequence $\boldsymbol\sigma= \{ i_1,\dots,i_n \}$ is expressed
in terms of $n$ densities $f_{\boldsymbol{\bar\sigma}^{(k)}}(x,t), \; k=1,\dots,n,$ whose indices
$\boldsymbol{\bar\sigma}^{(k)} = \{ i_1,\dots,i_{k-1},-i_k,i_{k+1},\dots,i_n \} , \; k=1,\dots,n,$ differ from the index
$\boldsymbol\sigma= \{ i_1,\dots,i_{k-1},i_k,i_{k+1},\dots,i_n \}$ in a single element only. In other words, the equation for arbitrary density 
$f_{\boldsymbol\sigma}(x,t)$ in (\ref{syst3}) with index $\boldsymbol\sigma = \{ i_1,\dots,i_n \} $ links it only with those densities whose indices 
are located from $\boldsymbol\sigma$ at distance 1 in the Hamming metric.

\section{Governing Equation}

\numberwithin{equation}{section}

Let $\bold\Xi_n = \{ \boldsymbol\sigma_1, \dots, \boldsymbol\sigma_{2^n} \}$ denote the ordered set consisting of $2^n$ sequences, each being of length $n$ and having the form $\boldsymbol\sigma_k = \{ i_1^{(k)},i_2^{(k)}, \dots, i_n^{(k)} \}, \; i_j^{(k)}=\pm 1, \; j=1,\dots,n, \; k=1,\dots,2^n, \; n\ge 2$.
The order in $\bold\Xi_n$ may be arbitrary, but fixed. For our purposes it is convenient to choose and fix the lexicographical order of the sequences in $\bold\Xi_n = \left\{ \boldsymbol\sigma_k = \{ i_1^{(k)},i_2^{(k)},\dots,i_n^{(k)} \}, \; i_j^{(k)}=\pm 1, \; j=1,\dots,n, \; k=1,\dots,2^n \right\} $, that is, the order
$$\aligned
\boldsymbol\sigma_1 & = \{ -1,-1,\dots,-1,-1 \}, \\
\boldsymbol\sigma_2 & = \{ -1,-1,\dots,-1,+1 \}, \\
\boldsymbol\sigma_3 & = \{ -1,-1,\dots,+1,-1 \}, \\
\boldsymbol\sigma_4 & = \{ -1,-1,\dots,+1,+1 \}, \\
\dots & \dots\dots\dots\dots\dots\dots\dots\\
\boldsymbol\sigma_{2^n} & = \{ +1,+1,\dots,+1,+1 \}.
\endaligned$$
Note that this lexicographical order is isomorphic to the binary one by the identification $-1\mapsto 0$ and $+1\mapsto 1$, however, for the sake of visuality, we keep
just the lexicographical order.

Let $\rho(\cdot,\cdot) : \bold\Xi_n\times\bold\Xi_n \to \{ 0,1,\dots,n \}$ be the Hamming metric. For arbitrary element $\boldsymbol\sigma_k\in\bold\Xi_n, \; k=1,\dots,2^n$, define a subset $\bold M_k\subset\bold\Xi_n$ by the formula:
$$\bold M_k = \{ \boldsymbol\sigma_s\in\bold\Xi_n : \rho(\boldsymbol\sigma_s,\boldsymbol\sigma_k) = 1  \}, \qquad k=1,\dots,2^n.$$
Identifying the notations $f_k(x,t)\equiv f_{\boldsymbol\sigma_k}(x,t), \; \bold c_k\equiv\bold c_{\boldsymbol\sigma_k}, \; k=1,\dots,2^n$, system (\ref{syst3}) can be represented in the following ordered form:
\begin{equation}\label{eq1}
\frac{\partial f_k(x,t)}{\partial t} =  - \bold c_k \frac{\partial f_k(x,t)}{\partial x} - \boldsymbol\lambda f_k(x,t)
+ \sum_{ \{ m \; : \; \boldsymbol\sigma_m\in\bold M_k \} } \lambda_m \; f_m(x,t) , \qquad k=1,\dots,2^n.
\end{equation}
The main subject of our interest is the sum of functions (\ref{syst2})
\begin{equation}\label{eq2}
p(x,t) = \sum_{k=1}^{2^n} f_k(x,t) ,
\end{equation}
which is the transition probability density of the process $L(t)$ defined by (\ref{syst1}).

Introduce into consideration the column-vector of dimension $2^n$
$$\bold f = \bold f(x,t) = \left( f_1(x,t), f_2(x,t), \dots, f_{2^n}(x,t) \right)^T$$
and the diagonal $(2^n\times 2^n)$-matrix differential operator
\begin{equation}\label{eq3}
\frak D_n = \text{diag} \{ A_k, \; k=1,\dots,2^n \} ,
\end{equation}
where $A_k, \;  k=1,\dots,2^n,$ are the differential operators
$$A_k = \frac{\partial}{\partial t} + \bold c_k \frac{\partial}{\partial x}, \qquad  k=1,\dots,2^n.$$
Define the scalar $(2^n\times 2^n)$-matrix $\bold\Lambda_n = \Vert\xi_{sm}\Vert, \; s,m=1,\dots,2^n,$ with the elements
\begin{equation}\label{eq4}
\xi_{sm}= \left\{ \aligned
\boldsymbol\lambda , \qquad & \text{if} \quad s=m,\\
-\lambda_k, \qquad & \text{if} \quad \rho(\boldsymbol\sigma_s,\boldsymbol\sigma_m)=1 \; \text{and} \; i_k^{(s)}\neq i_k^{(m)},\\
0, \qquad & \text{otherwise,} \endaligned \right. \qquad s,m=1,\dots,2^n.
\end{equation}
In other words, the matrix $\bold\Lambda_n$ has the following structure. All the diagonal elements are equal to $\boldsymbol\lambda$.
At the intersection of the $s$-th row and the $m$-th column (corresponding to the sequences $\boldsymbol\sigma_s = \{ i_1^{(s)},i_2^{(s)},\dots,i_n^{(s)} \}$
and $\boldsymbol\sigma_m = \{ i_1^{(m)},i_2^{(m)},\dots,i_n^{(m)} \}$, such that the Hamming metric between them is 1), the element $-\lambda_k$ is located, where
$k$ is the position number of the non-coinciding elements in these sequences $\boldsymbol\sigma_s$ and $\boldsymbol\sigma_m$. Note that, since the Hamming metric between these sequences is 1, such position number $k$ is unique. All other elements of the matrix are zeros. From this definition it follows that each row or column of matrix $\bold\Lambda_n$ contains $(n+1)$ non-zero elements and $2^n-(n+1)$ zeros. The sum of all the elements of every row or column of matrix $\bold\Lambda_n$ is zero 
in view of the definition of $\boldsymbol\lambda$ given by (\ref{lam}).

In these notations the system (\ref{eq1}) can be represented in the matrix form
\begin{equation}\label{eq5}
\left[ \frak D_n + \bold\Lambda_n \right] \bold f = \bold 0,
\end{equation}
where $\bold 0 = (0,0,\dots,0)$ is the zero vector of dimension $2^n$.

\bigskip

{\bf Theorem 2.} {\it The transition probability density $p(x,t)$ of the process $L(t)$ given by} (\ref{eq2}) {\it satisfies the following hyperbolic
partial differential equation of order $2^n$ with constant coefficients}
\begin{equation}\label{eq6}
\left\{ \text{Det} \left[ \frak D_n + \bold\Lambda_n \right]\right\} p(x,t) = 0,
\end{equation}
{\it where} $\text{Det} \left[ \frak D_n + \bold\Lambda_n \right]$ {\it is the determinant of the matrix differential operator}
$\left[ \frak D_n + \bold\Lambda_n \right] .$

\vskip 0.2cm

{\it Proof.} The proof immediately emerges by applying the Determinant Theorem \cite[the Theorem]{kol3}, \cite[Theorem 1]{kol4}
(see also \cite[Theorem 2]{kol1}) to system (\ref{eq5}). According to this Determinant Theorem, in order to extract the governing equation from a system of first-order differential equations with commuting differential operators, one only needs to compute the determinant of this system whose elements are the commuting differential operators (moreover, this theorem is also true in the case when the determinant is replaced by a polylinear form defined on an arbitrary commutative ring over the field of complex numbers). Since the differential operators $A_k, \; k=1,\dots,2^n,$ commute with each other, this Determinant Theorem is applicable to our case.

The hyperbolicity of equation (\ref{eq6}) follows from the hyperbolicity of system (\ref{syst3}) (or (\ref{eq5})). The theorem is proved. $\square$

\bigskip

{\it Remark 2.} Derivation of a general analytical formula for the determinant $\text{Det} \left[ \frak D_n + \bold\Lambda_n \right]$ is a fairly difficult algebraic  problem which lies apart of the purposes of this article. Nevertheless, this problem can considerably be simplified, if we notice that, from the form of system (\ref{syst3}), it follows that matrix $\frak D_n + \bold\Lambda_n$ has the block structure
\begin{equation}\label{eq7}
\frak D_n + \bold\Lambda_n =
\begin{pmatrix} \frak D_{n-1}^{(1)} + \bold\Lambda_{n-1} & \bold E_{n-1} \\
                \bold E_{n-1} & \frak D_{n-1}^{(2)} + \bold\Lambda_{n-1}
\end{pmatrix}  ,
\end{equation}
where the blocks in (\ref{eq7}) are composed of the following $(2^{n-1}\times 2^{n-1})$-matrices:
\begin{equation}\label{eq8}
\frak D_{n-1}^{(1)} = \text{diag} \{ A_k, \; k=1,\dots,2^{n-1} \} , \qquad \frak D_{n-1}^{(2)} = \text{diag} \{ A_k, \; k=2^{n-1}+1,\dots,2^n \} ,
\end{equation}
the $(2^{n-1}\times 2^{n-1})$-matrix $\bold\Lambda_{n-1}$ is defined similarly (\ref{eq4}) (but taking into account its dimension), and
$\bold E_{n-1} = -\lambda_1 E_{n-1} ,$ where $E_{n-1}$ is the unit matrix of dimension $(2^{n-1}\times 2^{n-1})$. Since the matrix $\bold E_{n-1}$ commute with $[\frak D_{n-1}^{(2)}+\bold\Lambda_{n-1}]$ (and, of course, with $[\frak D_{n-1}^{(1)}+\bold\Lambda_{n-1}]$ ), then applying the well-known Schur's formulas for the even-order determinants of block matrices to (\ref{eq7}), we obtain:
\begin{equation}\label{eq9}
\text{Det} \left[ \frak D_n + \bold\Lambda_n \right] = \text{Det} \left[ (\frak D_{n-1}^{(1)}+\bold\Lambda_{n-1}) (\frak D_{n-1}^{(2)}+\bold\Lambda_{n-1}) - \lambda_1^2 E_{n-1} \right] .
\end{equation}
Formula (\ref{eq9}) reduces computation of a determinant of dimension $(2^n\times 2^n)$ to the computation of a determinant of dimension $(2^{n-1}\times 2^{n-1})$. It is also useful to note that, in view of definition (\ref{eq4}), for arbitrary $n\ge 2$, the relation $\frak D_n + \bold\Lambda_n = \left[ \frak D_n + \bold\Lambda_n \right]^T$ holds. In other words, the matrix $\frak D_n + \bold\Lambda_n $ coincides with its transposed matrix.

Clearly, this approach can be extended recurrently to the determinants of lower dimensions with final obtaining an explicit (but complicated) formula for the determinant $\text{Det} \left[ \frak D_n + \bold\Lambda_n \right]$, however this is not our purpose here.

\bigskip

{\it Remark 3.} In the particular case $n=2$, we obtain the $(4\times 4)$-determinant
\begin{equation}\label{eq10}
\text{Det} \left[ \frak D_2 + \bold\Lambda_2 \right] =
\begin{vmatrix} A_1+\boldsymbol\lambda & -\lambda_2 & \vdots & -\lambda_1 & 0 \\
                -\lambda_2 & A_2+\boldsymbol\lambda & \vdots & 0 & -\lambda_1 \\
                 \hdotsfor{5} \\
                -\lambda_1 & 0 & \vdots & A_3+\boldsymbol\lambda & -\lambda_2 \\
                0 & -\lambda_1 & \vdots & -\lambda_2 & A_4+\boldsymbol\lambda
\end{vmatrix}
\end{equation}
In the case $n=3$, the following $(8\times 8)$-determinant emerges:
\begin{equation}\label{eq11}
\text{Det} \left[ \frak D_3 + \bold\Lambda_3 \right] =
\begin{vmatrix} A_1+\boldsymbol\lambda & -\lambda_3 & -\lambda_2 & 0 & \vdots & -\lambda_1 & 0 & 0 & 0 \\
                -\lambda_3 & A_2+\boldsymbol\lambda & 0 & -\lambda_2 & \vdots & 0 & -\lambda_1 & 0 & 0 \\
                -\lambda_2 & 0 & A_3+\boldsymbol\lambda & -\lambda_3 & \vdots & 0 & 0 & -\lambda_1 & 0 \\
                0 & -\lambda_2 & -\lambda_3 & A_4+\boldsymbol\lambda & \vdots & 0 & 0 & 0 & -\lambda_1 \\
                 \hdotsfor{9} \\
                -\lambda_1 & 0 & 0 & 0 & \vdots & A_5+\boldsymbol\lambda & -\lambda_3 & -\lambda_2 & 0 \\
                0 & -\lambda_1 & 0 & 0 & \vdots & -\lambda_3 & A_6+\boldsymbol\lambda & 0 & -\lambda_2 \\
                0 & 0 & -\lambda_1 & 0 & \vdots & -\lambda_2 & 0 & A_7+\boldsymbol\lambda & -\lambda_3 \\
                0 & 0 & 0 & -\lambda_1 & \vdots & 0 & -\lambda_2 & -\lambda_3 & A_8+\boldsymbol\lambda \\
\end{vmatrix}
\end{equation}
From (\ref{eq10}) and (\ref{eq11}) we clearly see the nice block structure of matrix $\frak D_n + \bold\Lambda_n$, as it was noted in (\ref{eq7}). It is also seen that
the diagonal blocks of determinant (\ref{eq11}) are structurally similar to (\ref{eq10}). Such determinants can, therefore, be evaluated by applying the 
recurrent formula (\ref{eq9}). Note that if we take some other order of the sequences $\{ \boldsymbol\sigma_k, \; k=1,\dots, 2^n\}$, (this corresponds to some pairwise changes of rows and columns in matrix $\frak D_n + \bold\Lambda_n$), the determinant $\text{Det} \left[ \frak D_n + \bold\Lambda_n \right]$ keeps the same value.

\bigskip

{\it Remark 4.} To obtain the fundamental solution of partial differential equation (\ref{eq6}) we should solve it with the initial conditions
\begin{equation}\label{eq12}
p(x,t)\vert_{t=0} = \delta\left( x - \sum_{k=1}^n a_k x_k^0 \right) , \qquad \left. \frac{\partial^k p(x,t)}{\partial t^k} \right|_{t=0} =0 , \quad k=1,\dots,2^n-1,
\end{equation}
where $\delta(\cdot)$ is the Dirac delta-function. The first condition in (\ref{eq12}) expresses the obvious fact that, at the initial time moment $t=0$, 
the density of process $L(t)$ is entirely concentrated at the start point $\sum\limits_{k=1}^n a_k x_k^0$. 

To pose the initial-value problem for the transition density $p(x,t)$ of process $L(t)$ we need to find the respective initial conditions. To do this, we may use 
the known formulas (\ref{prop10}) and (\ref{prop11}) for the characteristic function of the telegraph process. Since the telegraph processes $X_k(t), \; k=1,\dots,n,$ 
are independent, then, in view of (\ref{prop10}) and (\ref{prop11}), the characteristic function of their linear form $L(t)$ is given by the formula: 
\begin{equation}\label{charL}
H_L(\alpha,t) = \exp\biggl( -\boldsymbol\lambda t + i\alpha \sum_{k=1}^n a_k x_k^0 \biggr) \; \prod_{k=1}^n \tilde H_k(a_k \alpha,t) , 
\qquad \alpha\in\Bbb R, \quad t\ge 0, 
\end{equation}
where
\begin{equation}\label{charH}
\aligned
\tilde H_k(\xi, t) & = \left[ \cosh\left( t\sqrt{\lambda_k^2-c_k^2\xi^2} \right) + \frac{\lambda_k}{\sqrt{\lambda_k^2-c_k^2\xi^2}} \; 
\sinh\left( t\sqrt{\lambda_k^2-c_k^2\xi^2} \right) \right] \bold 1_{\left\{\vert\xi\vert\le\frac{\lambda_k}{c_k}\right\}} \\
& + \left[ \cos\left( t\sqrt{c_k^2\xi^2 - \lambda_k^2} \right) + \frac{\lambda_k}{\sqrt{c_k^2\xi^2-\lambda_k^2}} \;
\sin\left( t\sqrt{c_k^2\xi^2-\lambda_k^2} \right) \right] \bold 1_{\left\{\vert\xi\vert > \frac{\lambda_k}{c_k}\right\}} .
\endaligned
\end{equation}
In particular, setting $t=0$ in (\ref{charL}) we get the formula 
$$H_L(\alpha,0) = \exp\biggl( i\alpha \sum_{k=1}^n a_k x_k^0 \biggr)$$
and its inverting yields the first initial condition in (\ref{eq12}). To obtain other initial conditions, we should differentiate (in $t$) characteristic function 
(\ref{charL}) the respective number of times, then inverting (in $\alpha$) the result of such differentiation and setting then $t=0$. 

\bigskip

{\it Remark 5.} From the hyperbolicity of equation (\ref{eq6}) and initial conditions (\ref{eq12}) (more precisely, from the first inital condition of (\ref{eq12})) it follows that the fundamental solution $f(x,t)$ of equation (\ref{eq6}) is a generalized function and, therefore, the differential operator $\text{Det} \left[ \frak D_n + \bold\Lambda_n \right]$ in (\ref{eq6}) is treated, for any fixed $t>0$, as the differential operator acting in the space of generalized functions $\mathscr{S'}$. The elements of $\mathscr{S'}$ are called the {\it tempered distributions}. Such interpretation becomes more visual if we note that solving the initial-value problem (\ref{eq6})-(\ref{eq12}) is equivalent to solving the inhomogeneous equation
\begin{equation}\label{eq13}
\left\{ \text{Det} \left[ \frak D_n + \bold\Lambda_n \right]\right\} f(x,t) = \delta(t) \; \delta\left( x - \sum_{k=1}^n a_k x_k^0 \right) ,
\end{equation}
where the generalized function on the right-hand side of (\ref{eq13}) represents the weighted sum of the instant point-like sources concentrated, at the initial time moment $t=0$,  at the point $\sum\limits_{k=1}^n a_k x_k^0$. In such form of writing the initial-value problem (\ref{eq13}), the operator $\text{Det} \left[ \frak D_n + \bold\Lambda_n \right] : \mathscr{S'} \to \mathscr{S'}$ is the differential operator acting from $\mathscr{S'}$ to itself. From this point of view, solving the differential equation (\ref{eq6}) with initial conditions (\ref{eq12}) means finding a generalized function $f(x,t)\in \mathscr{S'}$ such that the differential operator $\text{Det} \left[ \frak D_n + \bold\Lambda_n \right]$ transforms it into the generalized function $\delta(t) \delta\left( x - \sum\limits_{k=1}^n a_k x_k^0 \right) \in \mathscr{S'}$. Since the initial-value problem (\ref{eq6})-(\ref{eq12}) is well-posed (due to the hyperbolicity of equation (\ref{eq6})), such generalized function $f(x,t)$ exists and is unique in $\mathscr{S'}$ for any fixed $t>0$. Therefore, the fundamental solution $f(x,t)$ of the linear form $L(t)$ defined by (\ref{syst1}) is the Green's function of the initial-value problem (\ref{eq6})-(\ref{eq12}). 

The same concerns the initial-value problem for the transition probability density $p(x,t)$ and the respective initial conditions are determined as described in Remark 4 above. Such initial-value problem can also be represented in the form of a inhomogeneous partial differential equation similar to (\ref{eq13}), but with another generalized function on its right-hand side determined by the initial conditions for the transition density $p(x,t)$.

\section{Limit Theorem}

\numberwithin{equation}{section}

In this section we examine the limiting behaviour of the linear form $L(t)$ defined by (\ref{syst1}) when the parameters of the telegraph processes tend to infinity 
in such a way that the following Kac's scaling conditions fulfil:
\begin{equation}\label{KacCond}
\lambda_k\to +\infty, \qquad c_k\to +\infty, \qquad \frac{c_k^2}{\lambda_k} \to\varrho_k, \quad \varrho_k>0,  \qquad k=1,\dots,n. 
\end{equation}
It is well known that, under condition (\ref{KacCond}), each the telegraph process $X_k(t)$ weakly converges to the homogeneous Wiener process $W_k(t)$ with zero drift and  diffusion coefficient $\sigma_k^2=\varrho_k$ starting from the initial point $x_k^0\in\Bbb R$. Therefore, it is quite natural to expect that the process $L(t)$ converges to 
the linear form 
$$W(t) = \sum_{k=1}^n a_k W_k(t)$$ 
of the Wiener processes $W_k(t), \; k=1,\dots, n$. In the following theorem we prove this fact. 

\bigskip

{\bf Theorem 3.} {\it Under Kac's scaling conditions} (\ref{KacCond}) {\it the weak convergence $L(t) \Rightarrow W(t)$ takes place, where $W(t)$ 
is the homogeneous Wiener process with the expectation and diffusion coefficient given, respectively, by the formulas:} 
\begin{equation}\label{Lim1}
\Bbb E W(t) = \sum_{k=1}^n a_k x_k^0 , \qquad \sigma_W^2 = \sum_{k=1}^n \varrho_k a_k^2 .
\end{equation} 

\vskip 0.2cm

{\it Proof}. Consider the characteristic function $H_L(\alpha,t)$ of the linear form $L(t)$ given by (\ref{charL}):
$$H_L(\alpha,t) = \exp\biggl( -\boldsymbol\lambda t + i\alpha \sum_{k=1}^n a_k x_k^0 \biggr) \; \prod_{k=1}^n \tilde H_k(a_k \alpha,t) , 
\qquad \alpha\in\Bbb R, \quad t\ge 0,$$
where, remind, $\boldsymbol\lambda = \sum\limits_{k=1}^n \lambda_k$ and functions $\tilde H_k(a_k \alpha,t)$ are given by the formula (see (\ref{charH})):
$$\aligned
& \tilde H_k(a_k \alpha, t) \\ 
& = \left[ \cosh\left( t\sqrt{\lambda_k^2-c_k^2a_k^2\alpha^2} \right) + \frac{\lambda_k}{\sqrt{\lambda_k^2-c_k^2a_k^2\alpha^2}} \; 
\sinh\left( t\sqrt{\lambda_k^2-c_k^2a_k^2\alpha^2} \right) \right] \bold 1_{\left\{\vert a_k \alpha\vert\le\frac{\lambda_k}{c_k}\right\}} \\
& + \left[ \cos\left( t\sqrt{c_k^2a_k^2\alpha^2 - \lambda_k^2} \right) + \frac{\lambda_k}{\sqrt{c_k^2a_k^2\alpha^2-\lambda_k^2}} \;
\sin\left( t\sqrt{c_k^2a_k^2\alpha^2-\lambda_k^2} \right) \right] \bold 1_{\left\{\vert a_k \alpha\vert > \frac{\lambda_k}{c_k}\right\}} ,
\endaligned$$
$$k=1,\dots, n.$$
From conditions (\ref{KacCond}) it follows that 
$$\frac{c_k^m}{\lambda_k^{m-1}} \to \left\{ \aligned \varrho_k, \quad & \text{if} \; m=2, \\
                                                    0, \quad & \text{if} \; m\ge 3, \endaligned \right.$$
for arbitrary $k=1,\dots, n$. Then the following asymptotic (under conditions (\ref{KacCond})) formula holds: 
$$\aligned 
t \sqrt{\lambda_k^2-c_k^2a_k^2\alpha^2} & = \lambda_k t \sqrt{1-\frac{c_k^2}{\lambda_k^2} \; a_k^2\alpha^2} \\
& = \lambda_k t \left[ 1 - \frac{1}{2} \; \frac{c_k^2}{\lambda_k^2} \; a_k^2\alpha^2 - \frac{1\cdot 1}{2\cdot 4} \; \left( \frac{c_k^2}{\lambda_k^2} \; a_k^2\alpha^2 \right)^2 - \frac{1\cdot 1\cdot 3}{2\cdot 4\cdot 6} \; \left( \frac{c_k^2}{\lambda_k^2} \; a_k^2\alpha^2 \right)^3 - \dots \right] \\
& = \lambda_k t - \frac{1}{2} \; \frac{c_k^2}{\lambda_k} \; a_k^2\alpha^2 t - \frac{1}{8} \; \frac{c_k^4}{\lambda_k^3} \; a_k^4\alpha^4 t - \frac{1}{16} \;  \frac{c_k^6}{\lambda_k^5} \; a_k^6\alpha^6 t - \dots \\
& \sim \; \lambda_k t - \frac{\varrho_k a_k^2\alpha^2}{2} \; t . 
\endaligned$$
Therefore, taking into account that, under conditions (\ref{KacCond}), $(\lambda_k/c_k) \to +\infty, \; (c_k^2/\lambda_k^2) \to 0,$ for any $k=1,\dots,n,$ we arrive at the following asymptotic formulas: 
$$\bold 1_{\left\{\vert a_k \alpha\vert\le\frac{\lambda_k}{c_k}\right\}} \; \to 1, \qquad \bold 1_{\left\{\vert a_k \alpha\vert > \frac{\lambda_k}{c_k}\right\}} \; \to 0,$$
$$\cosh\left( t\sqrt{\lambda_k^2-c_k^2a_k^2\alpha^2} \right) \; \sim \; \cosh\left( \lambda_k t  - \frac{\varrho_k a_k^2\alpha^2}{2} \; t \right),$$
$$\frac{\lambda_k}{\sqrt{\lambda_k^2-c_k^2a_k^2\alpha^2}} \; 
\sinh\left( t\sqrt{\lambda_k^2-c_k^2a_k^2\alpha^2} \right) \; \sim \;  \sinh\left( \lambda_k t   - \frac{\varrho_k a_k^2\alpha^2}{2} \; t  \right). $$

Thus, 
$$\tilde H_k(a_k \alpha, t) \sim \cosh\left( \lambda_k t  - \frac{\varrho_k a_k^2\alpha^2}{2} \; t \right) + \sinh\left( \lambda_k t   - \frac{\varrho_k a_k^2\alpha^2}{2} \; t  \right) = \exp\left( \lambda_k t  - \frac{\varrho_k a_k^2\alpha^2}{2} \; t \right),$$
and the following asymptotic formula holds: 
$$\prod_{k=1}^n \tilde H_k(a_k \alpha,t) \; \sim \;  \exp\left[ \boldsymbol\lambda t - \frac{1}{2} \left( \sum_{k=1}^n \varrho_k a_k^2 \right) \alpha ^2 t \right] \qquad (\text{under conditions (\ref{KacCond})}).$$
Therefore, we finally obtain (under conditions (\ref{KacCond})) the convergence
$$\lim\limits_{\substack{c_k, \lambda_k\to +\infty\\ (c_k^2/\lambda_k)\to\varrho_k}} H_L(\alpha,t) = \exp\left[ i\alpha \sum_{k=1}^n a_k x_k^0 - 
\frac{1}{2} \left( \sum_{k=1}^n \varrho_k a_k^2 \right) \alpha ^2 t \right] , \qquad k=1,\dots, n,$$
and the function on the right-hand side of this limiting relation is the characteristic function of the homogeneous Wiener process with the expectation and diffusion coefficient given, respectively, by formulas (\ref{Lim1}).
From this convergence of characteristic functions, the weak convergence $L(t) \Rightarrow W(t)$ follows. The theorem is proved. $\square$

\bigskip 

{\it Remark 6.} One can prove a more strong result concerning the pointwise convergence of the distribution of the linear form $L(t)$ to the distribution of $W(t)$ 
( that is, convergence at every point of the support of $L(t)$), but this requires a much more complicated analysis based on the theory of Cauchy problems in the space of generalized functions.

\section{Sum and Difference of Two Telegraph Processes}

\numberwithin{equation}{section}

In this section we apply the results obtained above for studying the sum and difference 
\begin{equation}\label{two1}
S^{\pm}(t) = X_1(t) \pm  X_2(t)
\end{equation}
of two independent telegraph processes $X_1(t)$ and $X_2(t)$.

The sum of two independent telegraph processes on the real line $\Bbb R$, both with the same parameters $c_1=c_2=c, \; \lambda_1=\lambda_2=\lambda$, that simultaneously  start from the origin $0\in\Bbb R$, was thoroughly studied in \cite{kol1} and the explicit probability distribution of this sum was obtained. It was also shown that the shifted time derivative of the transition density satisfies the telegraph equation with doubled parameters $2c$ and $2\lambda$. A functional relation connecting the distributions of the difference of two independent telegraph processes with arbitrary parameters and of the Euclidean distance between them, was given in \cite[Remark 4.4]{kol2}.

The results of the previous sections enable us to consider the generalizations of these models and to study the behaviour of the sum and difference of two independent telegraph processes $X_1(t), X_2(t)$ that, at the initial time moment $t=0$, simultaneously start from two arbitrary initial points $x_1^0, x_2^0\in\Bbb R$ and are developing with arbitrary constant velocities $c_1$ and $c_2$, respectively. The motions are controlled by two independent Poisson processes of arbitrary rates $\lambda_1$ and $\lambda_2$, respectively, as described above.

The coefficients of linear form (\ref{two1}) are $a_1=1, \; a_2=1$ for the sum $S^+(t)$ and $a_1=1, \; a_2=-1$ for the difference $S^-(t)$, respectively. Therefore, according to (\ref{supp2}), the supports of the distributions of $S^{\pm}(t)$ are the intervals
\begin{equation}\label{ttwo1}
\text{supp} \; S^{\pm}(t) = [ (x_1^0 \pm x_2^0)-(c_1+c_2)t, \; (x_1^0 \pm x_2^0)+(c_1+c_2)t ] .
\end{equation}
The lexicographically-ordered set of sequences in this case is
$$\boldsymbol\sigma_1=(-1,-1), \quad \boldsymbol\sigma_2=(-1,+1), \quad \boldsymbol\sigma_3=(+1,-1), \quad \boldsymbol\sigma_4=(+1,+1) ,$$
and according to (\ref{sing1}), the support of the sum $S^+(t)$ has, therefore, the following singularity points:
\begin{equation}\label{ttwo2}
\aligned 
q_1^+ & = (x_1^0+x_2^0)-(c_1+c_2)t, \qquad (\text{terminal point of the support}),\\
q_2^+ & = (x_1^0+x_2^0)-(c_1-c_2)t, \qquad (\text{interior point of the support}),\\
q_3^+ & = (x_1^0+x_2^0)+(c_1-c_2)t, \qquad (\text{interior point of the support}),\\
q_4^+ & = (x_1^0+x_2^0)+(c_1+c_2)t, \qquad (\text{terminal point of the support}).
\endaligned
\end{equation}
By setting $x_1^0=x_2^0=0, \; c_1=c_2=c, \; \lambda_1=\lambda_2=\lambda$, we arrive at the model studied in \cite{kol1} with the support 
$\text{supp} \; S^+(t) = [ -2ct, \; 2ct ]$. In this case formulas (\ref{ttwo2}) produce the 
three singularity points, namely, $\pm 2ct$ (the terminal points of the support) and $0$ (the interior point of multiplicity 2). 

Similarly, the support of the difference $S^-(t)$ has the following singularity points:
\begin{equation}\label{dif3}
\aligned 
q_1^- & = (x_1^0-x_2^0)-(c_1-c_2)t, \qquad (\text{interior point of the support}),\\
q_2^- & = (x_1^0-x_2^0)-(c_1+c_2)t, \qquad (\text{terminal point of the support}),\\
q_3^- & = (x_1^0-x_2^0)+(c_1+c_2)t, \qquad (\text{terminal point of the support}),\\
q_4^- & = (x_1^0-x_2^0)+(c_1-c_2)t, \qquad (\text{interior point of the support}).
\endaligned
\end{equation}

Note that if both the processes $X_1(t)$ and $X_2(t)$ start from the same initial point $x_1^0=x_2^0=x^0\in\Bbb R$ and have the same speed $c_1=c_2=c,$ then the 
support of the difference $S^-(t)$ takes the form $\text{supp} \; S^-(t) = [ -2ct, \; 2ct ]$ with the three singularity points $0, \pm 2ct$ 
(the interior singularity point $0$ has multiplicity 2). We see that in this case difference $S^-(t)$ has the same support and the same singularity points 
like the sum $S^+(t)$ of two telegraph processes with the same speed $c_1=c_2=c$ that simultaneously start from the origin $0\in\Bbb R$.

In view of (\ref{sing2}),
$$\text{Pr} \left\{ S^{\pm}(t) = q_j^{\pm} \right\} = \frac{e^{-\boldsymbol\lambda t}}{4} , \qquad j=1,2,3,4,$$
where $\boldsymbol\lambda = \lambda_1+\lambda_2$.   

According to (\ref{coeffC}), for the sum $S^+(t)$ the coefficients $\bold c_k^+ \equiv \bold c_{\boldsymbol\sigma_k}^+, \; k=1,2,3,4,$ are:
$$\bold c_1^+ = -(c_1+c_2), \quad  \bold c_2^+ = -(c_1-c_2), \quad \bold c_3^+ = c_1-c_2, \quad \bold c_4^+ = c_1+c_2 .$$
Then operators $A_k^+, \; k=1,2,3,4,$ take the form:
\begin{equation}\label{two2}
\aligned
A_1^+ & = \frac{\partial}{\partial t} -(c_1+c_2) \frac{\partial}{\partial x}, \qquad  A_2^+ = \frac{\partial}{\partial t} -(c_1-c_2) \frac{\partial}{\partial x}, \\
A_3^+ & = \frac{\partial}{\partial t} +(c_1-c_2) \frac{\partial}{\partial x}, \qquad A_4^+ = \frac{\partial}{\partial t} +(c_1+c_2) \frac{\partial}{\partial x} .
\endaligned
\end{equation}

Similarly, for the difference $S^-(t)$ the coefficients $\bold c_k^- \equiv \bold c_{\boldsymbol\sigma_k}^-, \; k=1,2,3,4,$ are:
$$\bold c_1^- = -(c_1-c_2), \quad  \bold c_2^- = -(c_1+c_2), \quad \bold c_3^- = c_1+c_2, \quad \bold c_4^- = c_1-c_2 ,$$
and, therefore, the operators $A_k^-, \; k=1,2,3,4,$ become 
\begin{equation}\label{dif4}
\aligned
A_1^- & = \frac{\partial}{\partial t} -(c_1-c_2) \frac{\partial}{\partial x}, \qquad  A_2^- = \frac{\partial}{\partial t} -(c_1+c_2) \frac{\partial}{\partial x}, \\
A_3^- & = \frac{\partial}{\partial t} +(c_1+c_2) \frac{\partial}{\partial x}, \qquad A_4^- = \frac{\partial}{\partial t} +(c_1-c_2) \frac{\partial}{\partial x} .
\endaligned
\end{equation}

The initial-value problems for the transition densities of processes (\ref{two1}) are given by the following theorem.

\bigskip

{\bf Theorem 4.} {\it The transition probability densities $p^{\pm}(x,t)$ of processes} (\ref{two1}) {\it are the solutions of the initial-value problems}
\begin{equation}\label{two3}
\aligned
& \biggl\{ \left( \frac{\partial}{\partial t}+(\lambda_1+\lambda_2)  \right)^2 \left[ \frac{\partial^2}{\partial t^2} + 2(\lambda_1+\lambda_2) \frac{\partial}{\partial t} - 2(c_1^2+c_2^2) \frac{\partial^2}{\partial x^2} - (\lambda_1-\lambda_2)^2 \right] \\
& \qquad\qquad\qquad\qquad\qquad + \left[ (c_1^2-c_2^2) \frac{\partial^2}{\partial x^2} + (\lambda_1^2-\lambda_2^2) \right]^2 \biggr\} p^{\pm}(x,t) = 0 , 
\endaligned
\end{equation}

\begin{equation}\label{incond}
\aligned 
p^{\pm}(x,t)\vert_{t=0} & = \delta\bigl( x - (x_1^0 \pm x_2^0) \bigr) , \qquad \frac{\partial p^{\pm}(x,t)}{\partial t}\biggl\vert_{t=0} = 0, \\ 
\frac{\partial^2 p^{\pm}(x,t)}{\partial t^2}\biggl\vert_{t=0} & = (c_1^2+c_2^2) \; \delta'' \bigl( x - (x_1^0 \pm x_2^0) \bigr) , \\ 
\frac{\partial^3 p^{\pm}(x,t)}{\partial t^3}\biggl\vert_{t=0} & = -2(\lambda_1c_1^2+\lambda_2c_2^2) \; \delta'' \bigl( x - (x_1^0 \pm x_2^0) \bigr) ,
\endaligned 
\end{equation}
{\it where $\delta''(x)$ is the second generalized derivative of Dirac delta-function.}

{\it Since equation} (\ref{two3}) {\it is hyperbolic, then, for arbitrary $t>0$, the solutions $p^{\pm}(x,t)$ of initial-value problems} (\ref{two3})-(\ref{incond}) {\it exist and are unique in the class of generalized functions $\mathscr{S'}$.}

\vskip 0.2cm

{\it Proof.} To begin with, we establish initial conditions (\ref{incond}). According to (\ref{charL}), the characteristic functions of the processes $S^{\pm}(t)$ are  
$$H^{\pm}(\alpha,t) = \exp\bigl( -(\lambda_1+\lambda_2)t + i\alpha(x_1^0 \pm x_2^0 \bigr) \; \tilde H_1(\alpha,t) \tilde H_2(\alpha,t) , 
\qquad \alpha\in\Bbb R, \quad t\ge 0, $$
where functions $\tilde H_1(\alpha,t), \; \tilde H_2(\alpha,t)$ are given by (\ref{charH}). Differentiating $H^{\pm}=H^{\pm}(\alpha,t)$ in $t$, 
after some calculations we obtain:
$$\aligned 
H^{\pm}(\alpha,t)\vert_{t=0} & = \exp\bigl( i\alpha(x_1^0 \pm x_2^0) \bigr) , \qquad\qquad  \frac{\partial H^{\pm}}{\partial t}\biggl\vert_{t=0} = 0, \\
\frac{\partial^2 H^{\pm}}{\partial t^2}\biggl\vert_{t=0} & = -(c_1^2+c_2^2) \alpha^2 e^{i\alpha(x_1^0 \pm x_2^0)} , \qquad \frac{\partial^3 H^{\pm}}{\partial t^3}\biggl\vert_{t=0} = 2(\lambda_1c_1^2+\lambda_2c_2^2) \alpha^2 e^{i\alpha(x_1^0 \pm x_2^0)}.
\endaligned$$
Inverting these functions in $\alpha$ yields initial conditions (\ref{incond}).

Let us now derive the governing equation for the transition density of the sum $S^+(t)$. To simplify the notations, we 
identify operators $A_k \equiv A_k^+, \; k=1,2,3,4$, by omitting the upper index, bearing in mind, however, that we deal with the operators $A_k^+$ presented by
(\ref{two2}). Thus, according to Theorem 2, we should evaluate determinant (\ref{eq10}) with operators $A_k$ given by (\ref{two2}). To do this, we apply formula 
(\ref{eq9}) to determinant (\ref{eq10}). We have: 
\begin{equation}\label{det2}
\aligned
\text{Det} \left[ \frak D_2 + \bold\Lambda_2 \right] & =
\begin{vmatrix} A_1+\boldsymbol\lambda & -\lambda_2 & \vdots & -\lambda_1 & 0 \\
                -\lambda_2 & A_2+\boldsymbol\lambda & \vdots & 0 & -\lambda_1 \\
                 \hdotsfor{5} \\
                -\lambda_1 & 0 & \vdots & A_3+\boldsymbol\lambda & -\lambda_2 \\
                0 & -\lambda_1 & \vdots & -\lambda_2 & A_4+\boldsymbol\lambda
\end{vmatrix} \\
\\
& = \text{Det} \left[
\begin{pmatrix} A_1+\boldsymbol\lambda & -\lambda_2 \\
                \\
                -\lambda_2 & A_2+\boldsymbol\lambda
\end{pmatrix}
\begin{pmatrix} A_3+\boldsymbol\lambda & -\lambda_2 \\
                 \\
                -\lambda_2 & A_4+\boldsymbol\lambda
\end{pmatrix} - \lambda_1^2
\begin{pmatrix} 1 \;\; & \; 0 \\
                \\
                0 \;\; & \; 1
\end{pmatrix} \right]
\\
\\
& =
\begin{vmatrix} (A_1+\boldsymbol\lambda)(A_3+\boldsymbol\lambda) - (\lambda_1^2-\lambda_2^2) & -\lambda_2(A_1+A_4+2\boldsymbol\lambda) \\
\\
                 -\lambda_2(A_2+A_3+2\boldsymbol\lambda) & (A_2+\boldsymbol\lambda)(A_4+\boldsymbol\lambda) - (\lambda_1^2-\lambda_2^2)
\end{vmatrix} ,
\endaligned
\end{equation}
where, remind, $\boldsymbol\lambda=\lambda_1+\lambda_2$. In view of (\ref{two2}), 
$$A_1 + A_4 = 2 \frac{\partial}{\partial t}, \qquad A_2 + A_3 = 2 \frac{\partial}{\partial t}$$ and, therefore, we have: 

\begin{equation}\label{two4}
\aligned
& \text{Det} \left[ \frak D_2 + \bold\Lambda_2 \right] \\
& =
\begin{vmatrix} (A_1+\boldsymbol\lambda)(A_3+\boldsymbol\lambda) - (\lambda_1^2-\lambda_2^2) & -2\lambda_2(\frac{\partial}{\partial t}+\boldsymbol\lambda) \\
\\
                 -2\lambda_2(\frac{\partial}{\partial t}+\boldsymbol\lambda) & (A_2+\boldsymbol\lambda)(A_4+\boldsymbol\lambda) - (\lambda_1^2-\lambda_2^2)
\end{vmatrix} \\
& = \bigl[ (A_1+\boldsymbol\lambda)(A_3+\boldsymbol\lambda) - (\lambda_1^2-\lambda_2^2) \bigr] \bigl[ (A_2+\boldsymbol\lambda)(A_4+\boldsymbol\lambda) - (\lambda_1^2-\lambda_2^2) \bigr] - 4\lambda_2^2\left( \frac{\partial}{\partial t}+\boldsymbol\lambda \right)^2 \\
& = (A_1+\boldsymbol\lambda)(A_2+\boldsymbol\lambda)(A_3+\boldsymbol\lambda)(A_4+\boldsymbol\lambda) - (\lambda_1^2-\lambda_2^2) \left[ (A_1+\boldsymbol\lambda)(A_3+\boldsymbol\lambda) + (A_2+\boldsymbol\lambda)(A_4+\boldsymbol\lambda) \right] \\
& + (\lambda_1^2-\lambda_2^2)^2 - 4\lambda_2^2\left( \frac{\partial}{\partial t}+\boldsymbol\lambda \right)^2 .
\endaligned
\end{equation}
According to (\ref{two2}), we have
$$\aligned
& (A_1+\boldsymbol\lambda)(A_3+\boldsymbol\lambda) = \left[ \left( \frac{\partial}{\partial t}+\boldsymbol\lambda  \right)^2  -(c_1^2-c_2^2) \frac{\partial^2}{\partial x^2} \right] - 2c_2 \; \frac{\partial}{\partial x} \left( \frac{\partial}{\partial t}+\boldsymbol\lambda  \right) , \\
& (A_2+\boldsymbol\lambda)(A_4+\boldsymbol\lambda) = \left[ \left( \frac{\partial}{\partial t}+\boldsymbol\lambda  \right)^2  -(c_1^2-c_2^2) \frac{\partial^2}{\partial x^2} \right] + 2c_2 \; \frac{\partial}{\partial x} \left( \frac{\partial}{\partial t}+\boldsymbol\lambda  \right) .
\endaligned$$
Therefore,
\begin{equation}\label{two5}
\aligned
(A_1+\boldsymbol\lambda)(A_2+\boldsymbol\lambda)(A_3+\boldsymbol\lambda)(A_4+\boldsymbol\lambda) & = \left[ \left( \frac{\partial}{\partial t}+\boldsymbol\lambda  \right)^2  -(c_1^2-c_2^2) \frac{\partial^2}{\partial x^2} \right]^2 - 4c_2^2 \; \frac{\partial^2}{\partial x^2} \left( \frac{\partial}{\partial t}+\boldsymbol\lambda  \right)^2 \\
& = \left[ \left( \frac{\partial}{\partial t}+\boldsymbol\lambda  \right)^2  -(c_1^2+c_2^2) \frac{\partial^2}{\partial x^2} \right]^2 - 4c_1^2c_2^2 \; \frac{\partial^4}{\partial x^4} ,
\endaligned
\end{equation}
and
\begin{equation}\label{two6}
(A_1+\boldsymbol\lambda)(A_3+\boldsymbol\lambda) + (A_2+\boldsymbol\lambda)(A_4+\boldsymbol\lambda) = 2\left[ \left( \frac{\partial}{\partial t}+\boldsymbol\lambda  \right)^2  -(c_1^2-c_2^2) \frac{\partial^2}{\partial x^2} \right] .
\end{equation}
Substituting (\ref{two5}) and (\ref{two6}) into (\ref{two4}) we obtain
$$\aligned
\text{Det} \left[ \frak D_2 + \bold\Lambda_2 \right] & = \left[ \left( \frac{\partial}{\partial t}+\boldsymbol\lambda  \right)^2  -(c_1^2+c_2^2) \frac{\partial^2}{\partial x^2} \right]^2 - 4c_1^2c_2^2 \; \frac{\partial^4}{\partial x^4} \\
& \quad - 2(\lambda_1^2-\lambda_2^2) \left[ \left( \frac{\partial}{\partial t}+\boldsymbol\lambda  \right)^2  -(c_1^2-c_2^2) \frac{\partial^2}{\partial x^2} \right] + (\lambda_1^2-\lambda_2^2)^2 - 4\lambda_2^2\left( \frac{\partial}{\partial t}+\boldsymbol\lambda \right)^2 \\
& = \left[ \left( \frac{\partial}{\partial t}+\boldsymbol\lambda  \right)^2  -(c_1^2+c_2^2) \frac{\partial^2}{\partial x^2} \right]^2 - 4c_1^2c_2^2 \; \frac{\partial^4}{\partial x^4} \\
& \quad - 2(\lambda_1^2+\lambda_2^2) \left( \frac{\partial}{\partial t}+\boldsymbol\lambda  \right)^2 + 2(\lambda_1^2-\lambda_2^2) (c_1^2-c_2^2)
\frac{\partial^2}{\partial x^2} + (\lambda_1^2-\lambda_2^2)^2 \\
\endaligned$$

$$\aligned
& = \left( \frac{\partial}{\partial t}+\boldsymbol\lambda  \right)^2 \left[ \left( \frac{\partial}{\partial t}+\boldsymbol\lambda \right)^2 - 2(\lambda_1^2+\lambda_2^2) \right] + \left[ (c_1^2+c_2^2)^2 - 4c_1^2c_2^2 \right] \frac{\partial^4}{\partial x^4} \\
& \qquad - 2(c_1^2+c_2^2) \left( \frac{\partial}{\partial t}+\boldsymbol\lambda  \right)^2 \frac{\partial^2}{\partial x^2} + 2(\lambda_1^2-\lambda_2^2) (c_1^2-c_2^2)
\frac{\partial^2}{\partial x^2} + (\lambda_1^2-\lambda_2^2)^2 \\
& = \left( \frac{\partial}{\partial t}+\boldsymbol\lambda  \right)^2 \left[ \frac{\partial^2}{\partial t^2} + 2\boldsymbol\lambda \frac{\partial}{\partial t} - (\lambda_1-\lambda_2)^2 \right] + (c_1^2-c_2^2)^2 \frac{\partial^4}{\partial x^4} \\
& \qquad - 2(c_1^2+c_2^2) \left( \frac{\partial}{\partial t}+\boldsymbol\lambda  \right)^2 \frac{\partial^2}{\partial x^2} + 2(\lambda_1^2-\lambda_2^2) (c_1^2-c_2^2)
\frac{\partial^2}{\partial x^2} + (\lambda_1^2-\lambda_2^2)^2 \\
& = \left( \frac{\partial}{\partial t}+\boldsymbol\lambda  \right)^2 \left[ \frac{\partial^2}{\partial t^2} + 2\boldsymbol\lambda \frac{\partial}{\partial t} - 2(c_1^2+c_2^2) \frac{\partial^2}{\partial x^2} - (\lambda_1-\lambda_2)^2 \right] \\
& \qquad + \left[ (c_1^2-c_2^2) \frac{\partial^2}{\partial x^2} + (\lambda_1^2-\lambda_2^2) \right]^2 ,
\endaligned$$
proving equation (\ref{two3}) for the transition density $p^+(x,t)$ of the sum $S^+(t)$. 

Comparing now operators (\ref{two2}) and (\ref{dif4}), we see that $A_1^- = A_2^+, \; A_2^- = A_1^+, \; A_3^- = A_4^+, \; A_4^- = A_3^+$. Therefore, as is easy to see, 
determinant (\ref{eq10}) written for operators $A_k^-, \; k=1,2,3,4,$ takes the same value like determinant (\ref{det2}) for operators $A_k^+, \; k=1,2,3,4$. 
Thus, equation (\ref{two3}) is also valid for the transition density $p^-(x,t)$ of the difference $S^-(t)$. This completes the proof of the theorem. $\square$

\bigskip

{\it Remark 7.} By setting $c_1=c_2=c$ and $\lambda_1=\lambda_2=\lambda$ in (\ref{two3}) we arrive at the fourth-order equation 
$$\left( \frac{\partial}{\partial t} + 2\lambda \right)^2 \left( \frac{\partial^2}{\partial t^2} + 4\lambda \frac{\partial}{\partial t} - 
4c^2 \frac{\partial^2}{\partial x^2} \right) p^{\pm}(x,t) = 0$$
and this result is weaker in comparison with the third-order equation obtained in \cite[formula (4.3) therein]{kol1} for the transition density 
of the sum $S^+(t)$ of two independent telegraph processes, both with the same parameters $(c, \lambda)$. It is interesting to note that in the product of these operators, the second one represents the Goldstein-Kac telegraph operator with doubled parameters $(2c, 2\lambda)$. 

\bigskip

{\it Remark 8.} The form of equation (\ref{two3}) enables us to make some interesting probabilistic observations. We see that the first term of equation (\ref{two3}) in square brackets represents a telegraph-type operator (containing also the free term $ - (\lambda_1-\lambda_2)^2$) which is invariant with respect to parameters $\lambda_1, \lambda_2$ and $c_1, c_2$. In other words, if we change $\lambda_1$ for  $\lambda_2$ and inversely, and/or $c_1$ for $c_2$ and inversely, the first telegraph-type term of equation (\ref{two3}) preserves its form, while the second one can vary its interior signs. From equation (\ref{two3}) we can conclude that the sum and difference $S^{\pm}(t)$ of two independent telegraph processes are not the telegraph processes, however they still contain some telegraph-type component. To show that, we represent equation (\ref{two3}) in the following form: 

\begin{equation}\label{two7}
\aligned
& \biggl\{ \left( \frac{\partial}{\partial t}+(\lambda_1+\lambda_2)  \right)^2 \left[ \frac{\partial^2}{\partial t^2} + 2(\lambda_1+\lambda_2) \frac{\partial}{\partial t} - 2(c_1^2+c_2^2) \frac{\partial^2}{\partial x^2} \right] \\
& \quad - \left[ (\lambda_1-\lambda_2) \frac{\partial}{\partial t} - (c_1^2-c_2^2) \frac{\partial^2}{\partial x^2}\right] \left[ (\lambda_1-\lambda_2) \frac{\partial}{\partial t} + (c_1^2-c_2^2) \frac{\partial^2}{\partial x^2} + 2(\lambda_1^2-\lambda_2^2) \right] \biggr\} p^{\pm}(x,t) = 0 .
\endaligned
\end{equation}
We see that the first term of (\ref{two7}) contains exactly the telegraph operator quite similar to the classical Goldstein-Kac operator (\ref{prop3}) with the replacements $\lambda\mapsto\lambda_1+\lambda_2$ and $c^2\mapsto 2(c_1^2+c_2^2)$. The second term of (\ref{two7}) represents the product of two heat operators and this fact implies the presence of a Brownian-type component in the processes $S^{\pm}(t)$. Notice also that in this product the first operator in square brackets is exactly the standard heat operator (for  $\lambda_1>\lambda_2$ and $c_1>c_2$), while the second one represents an inverse-time heat operator (that is, with the inverse time replacement $t\mapsto -t$) containing also the free term $2(\lambda_1^2-\lambda_2^2)$.

\bigskip

{\it Remark 9.} Solving the initial-value problem (\ref{two3})-(\ref{incond}) and obtaining the transition densities $p^{\pm}(x,t)$ of the sum and difference $S^{\pm}(t)$ 
of two independent telegraph processes with arbitrary parameters is a fairly difficult problem that will be realized in the framework of another project.

\numberwithin{equation}{section}

\end{document}